\documentclass[preprint,12pt]{elsarticle}
\usepackage[latin1]{inputenc}
\usepackage{subfigure}
\usepackage{amsmath,amssymb,amsthm}

\DeclareMathAlphabet{\mathpzc}{OT1}{pzc}{m}{it}
\newtheorem{thm}{Theorem}
\newtheorem{lemma}[thm]{Lemma}
\newtheorem{cor}[thm]{Corollary}
\newtheorem{prop}[thm]{Proposition}
\newtheorem{defn}[thm]{Definition}
\newtheorem{rem}[thm]{Remark}
\newtheorem{example}[thm]{Example}
\newtheorem{counterexample}[thm]{Counter-example}

\begin{document}

\begin{frontmatter}
\title{M\"obius inversion formula for monoids with zero}

\author[adr1]{Laurent Poinsot\corref{me}} 
\cortext[me]{laurent.poinsot@lipn.univ-paris13.fr}
\author[adr1]{Gérard H. E. Duchamp}
\author[adr1]{Christophe Tollu}
\address[adr1]{LIPN - UMR 7030, 
CNRS - Universit\'e Paris 13, 93430 Villetaneuse,
France}
\begin{abstract}
The M\"obius inversion formula, introduced during the 19th century in number theory, was generalized to a wide class of monoids called locally finite such as the free partially commutative, plactic and hypoplactic monoids for instance. In this contribution are developed and used some topological and algebraic notions for monoids with zero, similar to ordinary objects such as the (total) algebra of a monoid, the augmentation ideal or the star operation on proper series. The main concern is to extend the study of the M\"obius function to some monoids with zero, \textit{i.e.}, with an absorbing element, in particular the so-called Rees quotients of locally finite monoids. Some relations between the M\"obius functions of a monoid and its Rees quotient are also  provided.  
\end{abstract}

\end{frontmatter}



\section{Introduction}\label{section-introduction}

The classic M\"obius inversion formula from number theory, introduced during the 19th century, states that, for any complex or real-valued functions $f,g$ defined on the positive integers $\mathbb{N}\setminus\{0\}$, the following assertions are equivalent:
\begin{itemize}
\item For all $n$, $g(n)=\displaystyle\sum_{d|n}f(d)$.
\item For all $n$, $f(n)=\displaystyle\sum_{d|n}\mu(n/d)f(d)$.
\end{itemize}
In both formulae the sums are extended over all positive divisors $d$ of $n$, and $\mu$ is the classical M\"obius function. This result actually uses the fact that $\mu$ and $\zeta$ are inverse one from the other with respect to the usual Dirichlet convolution, where $\zeta$ is the characteristic function of positive integers (see for instance~\cite{Apo76}). 

This classic version of the M\"obius inversion formula was generalized in different ways by different authors. P. Doubilet, G.-C. Rota, and R. P. Stanley proposed a  systematic treatment of this problem for locally finite posets in~\cite{DRS75,STA96}, while P. Cartier and D. Foata in~\cite{cartfoa} proved such a formula holds in a wide class of monoids called \emph{locally finite}~\cite{Eil74}, and the M\"obius function was even explicitly computed for some of them. This paper is a contribution to the study of the M\"obius inversion formula, still in the context of locally finite monoids but for the particular case of monoids with zero. For instance, let $M$ be the set $\{0,1,a,b,c,ab,ac,ba,bc,ca,cb,abc,acb,bac,bca,cab,cba\}$. It becomes a monoid with zero when equipped with concatenation of words without common letters, also called \emph{standard words}; the other products give $0$. Let $\zeta_0$ be the  characteristic function of $M_0=M\setminus\{0\}$. Then, $\zeta_0$ is invertible -- with respect to convolution -- in the algebra  $\mathbb{Z}_0[M]$ of all functions that annihilate the zero $0$ of $M$, which is, in a first approximation, the $\mathbb{Z}$-algebra of polynomials in the noncommutative variables $\{a,b,c\}$ with only standard words as monomials. Indeed,  $\zeta_0=1+\zeta_0^+$, where $1$ is the characteristic function of the singleton $\{1\}$ and since $\zeta_0^+$ has no constant term, as a noncommutative polynomial (that is $\zeta_0^+(1)=0$), $\zeta_0$ is invertible, with inverse $\mu_0=\sum_{n\geq 0}(-\zeta_0^{+})^n$. Due to the particular multiplication in $M$, the `` proper part '' $\zeta_0^+$ of $\zeta_0$ is actually nilpotent, and the previous summation stops after four steps. Therefore $\mu_0$ can be computed by hand, and we obtain $\mu_0=1-a-b-c$.

Rather surprisingly, $\mu_0$ -- interpreted as the M\"obius function of the monoid with zero $M$ -- is the same as the M\"obius function of the free noncommutative monoid $\{a,b,c\}^*$. Moreover such a phenomenon also appears for less tractable monoids with zero: for instance, let us consider a monoid similar to $M$ but on an infinite alphabet $X$: it is the set of all words on $X$ without multiple occurrences of any letter, and with product $\omega\times\omega'$ equal to the usual concatenation $\omega\omega'$ when each letter appears at most one time in the resulting word, and $0$ otherwise.  Contrary to $M$, this monoid is found infinite. Nevertheless we can prove its characteristic function to be invertible, and its inverse is still equal to the usual M\"obius function of the free monoid $X^*$. In this case, it is not as easy to compute because the corresponding `` proper part '' is no more nilpotent, and the sum of a series needs to be evaluated in some relevant topology.

The explanation of this general phenomenon is given in the present paper whose main concern is the development of an algebraic and topological toolbox for a systematic and rigorous treatment of the M\"obius inversion formula for locally finite monoids with zero.

\section{Monoids with zero}

A monoid with zero is an ordinary monoid with a two-sided absorbing element, called the \emph{zero}. Such structures obviously occur in ring theory (the multiplicative monoid of an associative ring with unit is a monoid with zero), but they are also used to solve some (co)homological problems~\cite{Nov08,NP09}, and mainly in the study of ideal extensions of  semigroups~\cite{AS70,Cli65,CP61}.  

These structures are defined as follows: let $M$ be an ordinary monoid (with $1_M$ as its identity element) such that $|M|\geq 2$. Then, $M$ is called a \emph{monoid with zero} if, and only if, there is a two-sided absorbing element $0_M$,  \textit{i.e.}, $x 0_M=0_M=0_M x$ for every $x\in M$, with $0_M\not=1_M$. The distinguished element $0_M$ is called the \emph{zero} of $M$ (uniqueness is obvious). If in addition $M$ is commutative, then $M$ is called a \emph{commutative monoid with zero}. In the sequel, for any monoid $M$ with zero $0_M$, $M_{0}$ stands for $M\setminus\{0_M\}$.   

\begin{example}\label{expl-monZ}
\begin{enumerate}
\item The set of all natural numbers $\mathbb{N}$ with the ordinary multiplication is a commutative monoid with zero;
\item The multiplicative monoid of any (associative) ring $R$ with a unit $1_R$ is a monoid with zero $0_R$;
\item\label{pt3} If $M$ is any usual monoid (with or without zero), then for every $0\not \in M$, $M^{0}=M\cup\{0\}$ is a monoid with zero $0$: $x0=0=0x$ for every $x\in M^{0}$ extending the operation of $M$. It is commutative if, and only if, the same holds for $M$. The transformation of $M$ into $M^{0}$ is called an \emph{adjunction of a zero}, and $M^{0}$ is a monoid with a \emph{(two-sided) adjoined zero}. Note that $M^{0}$ is obviously isomorphic with $M^{\mathsf{z}}$ for every $\mathsf{z}\not\in M$, where $\mathsf{z}$ plays the same role as $0$;
\item\label{pt5} The set $\aleph_0\cup\{\aleph_0\}$ of all cardinal numbers less or equal to $\aleph_0$ (that is, the closed initial segment $[0,\aleph_0]$), with the usual cardinal addition (recall that $\aleph_0=[0,\aleph_0[=\mathbb{N}$ and $n+\aleph_0=\aleph_0=\aleph_0+n$ for every $n\leq \aleph_0$) is a commutative monoid with $\aleph_0$ as zero;
\item\label{smallcat} Let $\mathsf{C}$ be a small category~\cite{MLan98}. Then its set of arrows $\mathcal{A}(\mathsf{C})$, together with adjoined zero $0$ and identity $1$, is a monoid with zero when arrows composition is extended using 
$f\circ g=0$ whenever $\mathsf{dom}(f)\not=\mathsf{codom}(g)$ for every $f\in\mathcal{A}(\mathsf{C})$, and $f\circ 1=f=1\circ f$, $f\circ 0=0=0\circ f$ for every $f \in \mathcal{A}(\mathsf{C})\cup\{0,1\}$. Now suppose that $P$ is a poset, and $\mathsf{Int}(P)$ is the set of its intervals $[x,y]=\{z\in P : x\leq z \leq y\}$ for all $x\leq y$ in $P$ (see~\cite{DRS75,STA96}). An interval $[x,y]$ may be seen as an arrow from $x$ to $y$, and a composition may be defined: $[x,z]\circ[z,y]=[x,y]$. It follows that $P$ turns to be a small category, and $\mathsf{Int}(P)\cup\{0,1\}$, where $0,1\not\in \mathsf{Int}(P)$ and $0\not=1$,  becomes a monoid with zero. Another specialization is possible: let $n\in\mathbb{N}\setminus\{0\}$ be fixed, and consider the set $I$ of all pairs $(i,j)$ of integers such that $1\leq i,j\leq n$. Any usual $n$-by-$n$ \emph{matrix unit} $E_{(i,j)}$ may be seen as an arrow from $i$ to $j$, and such arrows are composed by $E_{(i,k)}\circ E_{(k,j)}=E_{(i,j)}$. Then $I$ becomes a small category, and the set of all matrix units, with adjoined $0$ and $1$, may be interpreted as a monoid with zero which is also quite similar to A.~Connes's groupoids~\cite{CON94}. 
\end{enumerate}
\end{example}

A major class of monoids with zero, that deserves a short paragraph on its own, is given by the so-called Rees quotients (see~\cite{AS70,CP61,Gri95}). Let $M$ be a monoid and $I$ be a two-sided ideal  of $M$, that is $IM\subseteq I \supseteq MI$, which is proper ($I$ is proper if, and only if, $I\not=M$, or, in other terms, $1_M\not\in I$). A congruence $\theta_I$ on $M$ is defined as follows: $(x,y)\in \theta_I$ if, and only if, $x,y\in I$ or $x=y$. The equivalence class of $x \in M$ modulo $\theta_I$ is \begin{math}\left \{ \begin{array}{ll}
\{x\} & \mbox{if $x\not\in I$;}\\
I & \mbox{if $x\in I$.}\end{array}\right . \end{math}
Therefore $I$ plays the role of a zero in the quotient monoid $M/\theta_I$, in such a way that it is isomorphic with the monoid with zero $(M\setminus I)\cup\{0\}$, where $0\not\in M\setminus I$, and with operation
\begin{equation}
x\times y = \left \{ \begin{array}{ll}
xy & \mbox{for $xy\not\in I$}\\
0 & \mbox{for $xy\in I$}
\end{array}\right .
\end{equation}
for every $x,y\in M\setminus I$, and $x\times 0=0=0\times x$ for every $x \in (M\setminus I)\cup\{0\}$. This monoid, unique up to isomorphism (the choice of the adjoined zero), is called the \emph{Rees quotient of $M$ by $I$}, and  denoted $M/I$. In what follows, we identify the carrier sets of both isomorphic monoids $M/\theta_I$ and $(M\setminus I)\cup\{0\}$, and we use juxtaposition for products in $M/I$ and in $M$. 
\begin{rem}
The fact that $I$ is proper guarantees that $1_M\in M\setminus I$, and therefore $1_M\not=0$.
\end{rem}
\begin{example}
Let $X=\{a,b,c\}$ and $I=\{\omega\in X^*: \exists x\in X,\ \mbox{such\ that}\ |\omega|_x\geq 2\}$, where $|\omega|_x$ denotes the number of occurrences of the letter $x$ in the word $\omega$. Then $X^*/I$ is the monoid with zero $M$ described in the Introduction (see Section~\ref{section-introduction}). 
\end{example}

\section{Contracted monoid algebra}

\emph{{\bf{Convention:}} In the present paper, a ring is assumed to be associative, commutative and with a unit $1_R$; the zero of a ring is denoted by $0_R$. An $R$-algebra $A$ is assumed to be associative (but non necessarily commutative) and has a unit $1$. Its zero is denoted by $0$.}

The main objective of this section is to recall the relevant version of the monoid algebra  of a monoid with zero over some given ring: in brief, the zeros of the monoid and the ring are identified.   
Let $R$ be a ring, and $X$ be any set. The \emph{support} of $f\in R^X$ is the set $\{x\in X : f(x)\not=0_R\}$. Now let $M$ be a monoid with zero $0_M$. Let us consider the usual monoid algebra $R[M]$ of $M$ over $R$, which is, as an $R$-module, the set $R^{(M)}$ of all maps from $M$ to $R$ with finite support, endowed with the usual Cauchy product~\cite{BouAlg}. By \emph{contracted monoid algebra} of $M$ over $R$ (see~\cite{CP61,Okn90}), we mean the factor algebra $R_0[M]=R[M]/R0_M$, where $R0_M$ is the two-sided ideal $R[(0_M)]=\{\alpha 0_M : \alpha \in R\}$. Thus, $R_0[M]$ may be identified with the set of all finite sums $\sum_{x\in M_0}\alpha_x x$, subject to the multiplication table given by the rule
\begin{equation}\label{constantes-de-structure-alg-contractee}
x\times y = \left \{
\begin{array}{lll}
xy & \mbox{if} & xy\not=0_M\ ,\\
0 & \mbox{if} & xy=0_M
\end{array}
\right .
\end{equation}
defined on basis $M_0$ (formula~(\ref{constantes-de-structure-alg-contractee}) gives the constants of structure, see~\cite{BouAlg}, of the algebra $R_0[M]$). In what follows we use juxtaposition rather than `` $\times$ '' for the products. From the definition, it follows directly that for any ordinary monoid $M$, $R_0[M^0]\cong R[M^0]/R0\cong R[M]$. This fact is extended to the Rees quotients as follows.
\begin{lemma}\label{1stlem}\cite{CP61,Okn90}
Let $M$ be a monoid and $I$ be a proper two-sided ideal of $M$. Then $R_0[M/I]\cong R[M]/R[I]$. (Note that $R[I]$ is the semigroup algebra of the ideal $I$.)
\end{lemma}
\begin{example}
Let $X$ be any non empty set and $n\in\mathbb{N} \setminus\{0,1\}$. Let $I$ be the proper ideal of $X^*$ of all words $\omega$ of length $|\omega|\geq n$. Then $R_0[X^*/I]$ consists in noncommutative polynomials truncated at length $n$.  
\end{example}
The notion of contracted monoid algebra is sufficient to treat the problem of the M\"obius formula for finite and locally finite (see Section~\ref{locallyfinitemonoidswithzero}) monoids with zero. Nevertheless infinite monoids with zero also occur, and formal series must be considered in those cases.

\section{Total contracted algebra of a finite decomposition monoid with zero}

Let $R$ be a ring, and $M$ be a usual monoid. The set of all functions $R^M$ has a natural structure of  $R$-module. By abuse of notation\footnote{\label{unenote}When $R^M$ is endowed with the topology of simple convergence, $R$ being discrete, the family $(f(x)x)_{x\in M}$ is summable, and $f=\sum_{x\in M}f(x)x$.}, any function $f\in R^M$ may be denoted by $\sum_{x\in M}\langle f,x \rangle x$, where\footnote{The notation `` $\langle f,x\rangle$ '' is commonly referred to as a `` Dirac bracket ''. It was successfully used by Sch\"utzenberger to develop his theory of automata~\cite{BR88}.} $\langle f,x\rangle=f(x)=\pi_x(f)$ ($\pi_x$ is the projection onto $Rx$). The carrier structure of the algebra $R[M]$ of the monoid $M$ is then seen as a submodule of $R^M$. Now taking $M$ to be a monoid with zero, we can also construct $R[M]$, however we would like to identify $0_M$ with $0$ of $R^M$ in the same way as  $R_0[M]$. Let us consider the set 
$R0_M=\{f\in R^M : \forall x\not=0_M,\ \langle f,x\rangle=0_R\}$,  \textit{i.e.}, $R0_M$ is the cyclic submodule generated by $0_M$. Then the quotient module $R^M/R0_M$ may be identified with the $R$-module $R^{M_0}$ of all `` infinite '' sums\footnote{As in the previous note~\ref{unenote}, it can be easily proved that such sums are actually the sums of summable series in the product topology on $R^M/R0_M$, with $R$ discrete.} $\sum_{x\in M_0}\langle f,x\rangle x$, or more likely the space of all functions from $M_0$ to $R$, \textit{i.e.}, $R^{M_0}=\{f\in R^M : f(0_M)=0_R\}$. This quotient module is the completion $\widehat{R_0[M]}$ of the topological module $R_0[M]$ equipped with the product topology ($R$ is given the discrete topology), also called `` topology of simple convergence '' or `` finite topology ''. It should be noticed that the quotient topology of $R^{M_0}$ induced by $R^M$ is equivalent to its product topology. 

Recall that an ordinary semigroup  (resp. monoid) $M$ is said to be a \emph{finite decomposition semigroup} (resp. \emph{finite decomposition monoid}), or to have the \emph{finite decomposition property}, if, and only if, it satisfies the following condition 
\begin{equation}\label{D-condition}
\forall x\in M,\ |\{(y,z)\in M\times M : yz=x\}|<+\infty\ .
\end{equation}
This condition is called the \emph{(D) condition} in~\cite{BouAlg}. If (\ref{D-condition}) holds, then $R^M$ can be equipped with the usual Cauchy or convolution product: therefore the $R$-algebra $R[[M]]$ of all formal power series over $M$ with coefficients in $R$ is obtained, which is also called the \emph{total algebra of the semigroup} (resp. \emph{monoid}) $M$ over $R$. This notion is now adapted to the case of monoids with zero. 
\begin{defn}
A monoid $M$ with zero $0_M$ is said to be a \emph{finite decomposition monoid with zero} if, and only if, it satisfies the following condition
\begin{equation}
\forall x\in M_0=M\setminus\{0_M\},\ |\{(y,z)\in M\times M : yz=x\}|<+\infty\ .
\end{equation}
\end{defn}
\begin{example}\label{locallyfiniteposet2finitedecompmonwithzero}
Let $P$ be a locally finite poset (\cite{DRS75,STA96}), \textit{i.e.}, such that every interval $[x,y]\in\mathsf{Int}(P)$ is finite. Then the monoid $\mathsf{Int}(P)\cup\{0,1\}$ of example~\ref{expl-monZ}.\ref{smallcat} is a finite decomposition monoid with zero. 
\end{example}
Some obvious results are given below without proofs.
\begin{lemma}\label{fd2fd0}
\begin{enumerate}
\item Let $M$ be a monoid with zero which has the finite decomposition property as an ordinary monoid. Then $M$ is finite.
\item Suppose that $M$ is a finite decomposition monoid. Then $M^0$ is a finite decomposition monoid with zero.
\item Suppose that $M$ is a finite decomposition monoid and $I$ is a two-sided proper ideal of $M$. Then the Rees quotient monoid $M/I$ is a finite decomposition monoid with zero. 
\end{enumerate}
\end{lemma}
Let us suppose that $M$ is a finite decomposition monoid with zero. Let $f,g\in R^M/R0_M$. Then we can define the corresponding Cauchy product:
\begin{equation}\label{Cauchyproductseries}
fg=\displaystyle\sum_{x\in M_0}\left (\sum_{yz=x}\langle f,y\rangle\langle g,z\rangle\right )x\ .
\end{equation}
The algebra $R^M/R0_M$ is then denoted $R_0[[M]]$ and called the \emph{total contracted algebra of the monoid $M$ over $R$}. The $R$-module $R_0[[M]]$ is the completion of $R_0[M]$ and because the Cauchy product of `` formal series '' in $R_0[[M]]$ is the continuous extension of its polynomial version in $R_0[M]$ (this product is separately continuous and continuous at zero~\cite{BouGenTop}), the following lemma holds. 
\begin{lemma}\label{completionalgebra4fdMw0}
Let $M$ be a finite decomposition monoid with zero. 
Then $R_0[[M]]$ is the completion of the contracted algebra $R_0[M]$, and, in particular, $R_0[[M]]$ is a topological algebra. 
\end{lemma}
Let $M$ be an ordinary monoid and $I$ be a two-sided proper ideal of $M$. Then the $R$-module $R^{M/I}/R0$ is isomorphic to the set of all formal infinite $R$-linear combinations $\sum_{x\not\in I}\langle f,x\rangle x$, where $f\in R^M$.  Now suppose that $M$ is a finite decomposition monoid. According to Lemma~\ref{fd2fd0}, $M/I$ is a finite decomposition monoid with zero. We can define both total algebras $R[[M]]$ and $R_0[[M/I]]$, with respectively  $R^M$ and $R^{M/I}/R0$ as carrier sets. The product on $R^{M/I}/R0$ is therefore given by 
\begin{equation}
\displaystyle\left ( \sum_{x\not\in I}\langle f,x\rangle x\right)\left ( \sum_{x\not\in I} \langle g,x\rangle x\right)=
\sum_{x\not\in I}\left ( \sum_{yz=x}\langle f,y\rangle \langle g,z\rangle\right )x\ .
\end{equation}
Let define 
\begin{equation}\label{epimorphisme}
\begin{array}{llll}
\Phi :& R[[M]]&\rightarrow &R_0[[M/I]]\\
&\displaystyle\sum_{x\in M}\langle f,x\rangle x &
\mapsto & \displaystyle \sum_{x\not\in I}\langle f,x\rangle x\ .
\end{array}
\end{equation}
Then $\Phi$ is an $R$-algebra homomorphism, which is onto and obviously continuous (for the topologies of simple convergence). Moreover $\ker(\Phi)=R[[I]]$, then $R_0[[M/I]]\cong R[[M]]/R[[I]]$. According to lemma~\ref{completionalgebra4fdMw0},  
$R_0[[M/I]]$ is complete (as an $R$-algebra) for the product topology. In summary we obtain:
\begin{prop}
Let $M$ be a finite decomposition monoid and $I$ be a proper two-sided ideal of $M$. Then, 
\begin{equation}
\begin{array}{lll}
R_0[[M/I]]&\cong& \widehat{R_0[M/I]}\\
&\cong& R[[M]]/R[[I]]\ .
\end{array}
\end{equation}
\end{prop}

\section{Locally finite monoids with zero}\label{locallyfinitemonoidswithzero}

In order to study the M\"obius inversion formula for monoids with zero, we need to characterize invertible series in the total contracted algebra. This can be done by exploiting a star operation on series without constant terms (\textit{i.e.}, for which $\langle f,1_M\rangle=0$). This star operation is easily defined when a topology on the algebra of series is given by some  filtration which generalizes the ordinary valuation. A particular class of monoids with zero satisfies this requirement. First we recall some classic results, and then we mimic them in the context of monoids with zero. 

A locally finite monoid $M$~\cite{cartfoa,Eil74} is a monoid such that
\begin{equation}
\forall x\in M,\  |\{ (n,x_1,\cdots, x_n) : x = x_1 \cdots x_n,\ x_i\not=1_M \}| < +\infty\ .
\end{equation}
For instance, any free partially commutative monoid~\cite{cartfoa,DK} is locally finite. 
A locally finite monoid is obviously a finite decomposition monoid, but the
converse is false since every non trivial finite group has the finite decomposition property, but is not locally finite because it has torsion. Furthermore, in a locally finite monoid, $xy=1_M\ \Rightarrow\ x=y=1_M$, or in other terms, $M\setminus\{1_M\}$ is a semigroup (and actually a locally finite semigroup in a natural sense), or, equivalently, the only invertible element of $M$ is the identity (such monoids are sometimes called \emph{conical}~\cite{COH03}). 
\begin{rem}
In~\cite{Bro68,She65} the authors -- L.N. Shevrin and T.C. Brown -- used another notion: they called \emph{locally finite} any semigroup in which every finitely-generated sub-semigroup is finite. This concept is really different and not comparable from the one used in this paper which follows~\cite{Eil74}.
\end{rem}
When M is locally finite, the $R$-algebra $R[[M]]$ may be equipped with a \emph{star
operation} defined for every \emph{proper series} $f$ (\textit{i.e.} such that $\langle f,1_M\rangle = 0_R$) by 
$f^{*} = \sum_{x\in M} \left( \sum_{n\geq 0}\sum_{x_1\cdots x_n=x}\langle f,x_1\rangle \cdots \langle f,x_n\rangle \right )x$ 
(\textit{i.e.} by $f^* = \sum_{n\in \mathbb{N}}f^n$). It follows that the \emph{augmentation ideal} $\mathfrak{M}=\{f\in R[[M]] : f\ \mbox{is\ proper}\}$, kernel of the usual augmentation map $\epsilon(f)=\langle f,1_M\rangle$ for every $f \in R[[M]]$, has the property that $1+\mathfrak{M}$ is a group (under multiplication; the inverse of $1-f\in 1+\mathfrak{M}$, when $f$ is proper, is precisely $f^*$), called the \emph{Magnus group} (see~\cite{DK} for instance). 
For this kind of monoids, we can define a natural notion of \emph{order function}. Let $x\in M$, then  $\omega_M (x)=\max\{n\in \mathbb{N} : 
\exists x_1,\cdots,x_n\in M\setminus\{1_M\},\ x=x_1\cdots x_n\}$. For instance if $M$ is a free partially commutative monoid $M(X,C)$, then $\omega_{M}(w)$ is the length $|w'|$  of $w'\in X^*$ of any element $w'$ in the class $w$. 

Let us adapt this situation to the case of monoids with zero. In what follows, if $M$ is any monoid (ordinary or with zero), then $M^{+}=M\setminus\{1_M\}$. A \emph{locally finite monoid with zero} (see~\cite{Esi08} for a similar notion) is a monoid with zero $M$ such that
\begin{equation}
\forall x\in M_0,\ |\{ (n,x_1,\cdots, x_n) : x = x_1 \cdots x_n,\ x_i\not=1_M \}| < +\infty\ .
\end{equation}
A locally finite monoid with zero obviously is also a finite decomposition monoid with zero. As in the case of usual monoids, the converse is false. Besides, if $M$ is a locally finite monoid, and $I$ is a two-sided proper ideal of $M$, then the Rees quotient $M/I$ is a locally finite monoid with zero. 
\begin{example}
Let $M=X^*/I$. Then $\omega_{M/I}(w)=|w|$ for every $w\in X^*\setminus I$. 
\end{example}
\begin{counterexample}
The monoid with zero $\mathsf{Int}(P)\cup\{0,1\}$ of a non void locally finite poset is not a locally finite monoid with zero, since for every $x \in P$, $1\not=[x,x]=[x,x]\cdot[x,x]$ holds.   
\end{counterexample}
As in the classical case, we can define a natural notion of \emph{order function} in a locally finite monoid with zero: let $x\in M_0$, then $\omega_M (x)=\max\{n\in \mathbb{N} : 
\exists x_1,\cdots,x_n\in M^{+},\ x=x_1\cdots x_n\}$ (we use the notation `` $\omega (x)$ '' when no confusion arises). Therefore $\omega (x)=0$ if, and only if, $x=1_M$. Moreover for every $x\in M_0$, if $x=yz$, then $\omega(x)\geq \omega(y)+\omega(z)$. If $M$ is a locally finite monoid and $I$ is a two-sided proper ideal of $M$, then we already know that $M/I$ is a locally finite monoid with zero, and more precisely for every $x\in M\setminus I$, $\omega_{M/I}(x)=\omega_M(x)$. 

Now let, $f\in R_0[[M]]$ (the total contracted algebra exists because $M$ is a finite decomposition monoid with zero since it is a locally finite monoid with zero). We define an \emph{order function} or \emph{pseudo-valuation} (that extends the order function $\omega_M$ of $M$): $\omega(f)=\inf\{\omega_M(x) : x\in M_0,\ \langle f,x\rangle\not=0_R\}$, where the infimum is taken in $\mathbb{N}\cup\{+\infty\}$. In particular, $\omega(f)=+\infty$ if, and only if, $f=0$. The following holds:
\begin{enumerate}
\item $\omega(1)=0$;
\item $\omega(f+g)\geq \min\{\omega(f),\omega(g)\}$;
\item $\omega(fg)\geq \omega(f)+\omega(g)$.
\end{enumerate}
Let us introduce $\mathfrak{M}=\{f\in R_0[[M]] : \langle f,1_M\rangle = 0_R\}=\{f\in R_0[[M]] : \omega(f)\geq 1\}$. This set obviously is a two-sided ideal of 
$R_0[[M]]$, called -- as in the ordinary case -- the \emph{augmentation ideal}\footnote{It is the kernel of the character $\epsilon : R_0[[M]]\rightarrow R$ given by $\epsilon(f)=\langle f,1_M\rangle=\pi_{1_M}(f)$, for $f \in R_0[[M]]$.}. For each $n\in\mathbb{N}$, let 
$\mathfrak{M}_{\geq n}= \{f \in R_0[[M]] : \omega(f)\geq n\}$, in such a way that $\mathfrak{M}_{\geq 0}=R_0[[M]]$, and 
$\mathfrak{M}_{\geq 1}=\mathfrak{M}$. The following lemma holds trivially. 
\begin{lemma}\label{filtration}
For every $n$, $\mathfrak{M}_{\geq n}$ is a two-sided ideal of $R_0[[M]]$, and 
the sequence $(\mathfrak{M}_{\geq n})_n$ is an exhaustive and separated decreasing filtration on $R_0[[M]]$, \textit{i.e.}, 
$\mathfrak{M}_{\geq n+1}\subseteq \mathfrak{M}_{\geq n}$, $\displaystyle\bigcup_{n\geq 0}\mathfrak{M}_{\geq n}=R_0[[M]]$, and 
$\displaystyle\bigcap_{n\geq 0}\mathfrak{M}_{\geq n}=(0)$. 
\end{lemma}

According to Lemma~\ref{filtration}, $R_0[[M]]$ with the topology $\mathcal{F}$ defined by the filtration $(\mathfrak{M}_{\geq n})_n$ is an Hausdorff topological ring (note also that this topology is metrizable~\cite{BouAlgCom}), and even an Hausdorff topological $R$-algebra when $R$ is discrete. 
\begin{rem}\label{equivalencetopologies}
It can be proved that if for every $n\in \mathbb{N}$, $M(n)=\{x\in M_0 : \omega_M (x)=n\}$ is finite, then the topology of simple convergence and the topology $\mathcal{F}$ on $R_0[[M]]$ are equivalent. In all cases, the topology defined by the filtration is always finer than the product topology (in particular, each projection $\pi_x : R_0[[M]]\rightarrow R$ is continuous for the filtration), and it can be even strictly finer as it is shown in the following example. 
\end{rem}
\begin{example}
Let us consider a countable set $X=\{x_i\}_{i\in\mathbb{N}}$ (that is $x_i\not=x_j$ for every $i\not=j$). We consider $M$ as the monoid $X^*$ with some zero $0$ adjoined. It is obviously a locally finite monoid with zero but the number of elements of a given order is not finite (for instance the number of elements of order $1$ is $\aleph_0$). We denote by $|\omega|$ the usual length of a word in $X^*$. Now let us consider the sequence of series $f_n = \sum_{k=0}^n x_k \in R_0[M]\subset R_0[[M]]$ which converges to the sum $f=\sum_{k=0}^{\infty}x_k$ in $R_0[[M]]$ endowed with the product topology ($f$ is the characteristic function of the alphabet $X$). But this series does not converge in $R_0[[M]]$ with the topology defined by the filtration, because $\omega(f-f_n)=1$ for all $n$. Nevertheless $f$ belongs to $R_0[[M]]$ since it is the completion of $R_0[M]$ in the product topology. 
\end{example}
Without technical difficulties the lemma below is obtained. 
\begin{lemma}
The algebra $R_0[[M]]$ with the topology $\mathcal{F}$ is complete. 
\end{lemma}
\begin{rem}
Suppose that $M$ is a locally finite monoid (with or without zero) which is also finite, then there exists $N\in\mathbb{N}$ such that for every $n\geq N$, $\mathfrak{M}_{\geq n}=(0)$.  In this case, the topology defined by the filtration coincides with the discrete topology on $R[[M]]=R[M]$ (or $R_0[[M]]=R_0[M]$). So no topology is needed in this case as explained in Introduction (Section~\ref{section-introduction}).
\end{rem}

\section{Star operation and the M\"obius inversion formula}

In this section, $M$ is assumed to be a locally finite monoid with zero. 
\begin{lemma}\label{lemmeinversionseriepropre}
For every $f\in\mathfrak{M}$, $(1-f)$ is invertible and $(1-f)^{-1}=\sum_{n\geq 0}f^n$.
\end{lemma}
\begin{proof}
First of all, $\sum_{n=0}^{+\infty}f^n$ is convergent in $R_0[[M]]$ (in the topology defined by the filtration), and is even summable, because $\omega(f^n)\rightarrow\infty$ when  $n\to+\infty$ (see~\cite{BouAlgCom}). Now for every $N\in \mathbb{N}$, $(1-f)\sum_{n=0}^N=1-f^{N+1}\rightarrow 1$ when $N\to+\infty$. Since $\sum_{n\geq 0}f^n$ is summable, and $R_0[[M]]$ is a topological algebra, we obtain asymptotically $(1-f)\sum_{n\geq 0}f^n=1$.  
\end{proof}

According to Lemma~\ref{lemmeinversionseriepropre}, for every element $f \in\mathfrak{M}$, we can define, as in the ordinary case, the \emph{star operation} $f^{*}=\sum_{n\geq 0}f^n$. 
\begin{rem}
Suppose that $M$ is a locally finite monoid with zero which is also finite. Then for every $f\in \mathfrak{M}$, $f$ is nilpotent (since $(f^n)_{n\in\mathbb{N}}$ is summable in the discrete topology). So in this particular case, there is no need of topology to compute $f^*$, as the example given in the Introduction. 
\end{rem}
\begin{lemma}
The set $1+\mathfrak{M}$ is a group under multiplication.
\end{lemma}
\begin{proof}
It is sufficient to prove that $\langle f^*,1_M\rangle=1_R$ for every $f\in\mathfrak{M}$. For every $n>0$, $\langle f^n,1_M\rangle=0$. Since the projection $\pi_{1_M}$ is continuous, we have 
\begin{equation}\langle f^*,1_M\rangle=\langle \displaystyle 1+\sum_{n\geq 1}f^n,1_M\rangle=\langle 1,1_M\rangle+\sum_{n\geq 0}\langle f^n,1_M\rangle=1_R\ .
\end{equation} 
\end{proof}

If $M$ is an ordinary locally finite monoid, the \emph{characteristic series} of $M$ is define as the series $\zeta=\sum_{x\in X}x \in R[[M]]$. If $X\subseteq M$, then $\underline{X}=\sum_{x\in X}x$ is the \emph{characteristic series} of $X$. More generally, if $M$ is a locally finite monoid with zero, then we also define the \emph{characteristic series} of $M$ by $\zeta_0=\sum_{x\in M_0}x\in R_0[[M]]$, and if $X\subseteq M$, then its characteristic series is $\underline{X}_0=\sum_{x\in X_0}x$ where $X_0=X\setminus\{0_M\}$.  We are now in a position to state the M\"obius inversion formula in the setting of (locally finite) monoids with zero.
\begin{prop}[M\"obius inversion formula]
The characteristic series $\zeta_0$ is invertible. 
\end{prop}
\begin{proof}
It is sufficient to prove that $\zeta_0\in 1+\mathfrak{M}$, which is obviously the case since $\zeta_0=1+\zeta_0^{+}$, where $\zeta_0^{+}=\underline{M^+}_0=\sum_{x\in M_0\setminus\{1_M\}}x\in \mathfrak{M}$. 
\end{proof}

We now apply several of the previous results on Rees quotients. So let $M$ be a locally finite monoid and $I$ be a two-sided proper ideal of $M$ in such a way that $M/I$ is a locally finite monoid with zero. 
Let us denote by $\mathfrak{M}_I$ (resp. $\mathfrak{M}$) the augmentation ideal of $M/I$ (resp. $M$). Let $\Phi : R[[M]]\rightarrow R_0[[M/I]]$ be the $R$-algebra epimorphism defined in eq.~(\ref{epimorphisme}). We know that it is continuous when both $R[[M]]$ and $R_0[[M/I]]$ have their topology of simple convergence. It is also continuous for the topologies defined by the filtrations $(\mathfrak{M}_{\geq n})_n$ and $((\mathfrak{M}_{I})_{\geq n})_n$. Indeed, let $n\in\mathbb{N}$, then for every $f\in \mathfrak{M}_{\geq n}$, $\Phi(f)\in(\mathfrak{M}_{I})_{\geq n}$. It admits a section $s$ from $R_0[[M/I]]$ into $R[[M]]$ (so $\Phi(s(f))=f$ for every $f \in R_0[[M/I]]$) defined by \begin{math}\langle s(f),x\rangle = \left \{\begin{array}{ll}
\langle f,x\rangle & \mbox{if}\ x\not\in I\ ,\\
0_R & \mbox{otherwise}\ . 
\end{array}
\right .\end{math}
This map is easily seen as an $R$-module morphism but in general not a ring homomorphism.
\begin{lemma}\label{interlemma1}
Let $f\in 1+\mathfrak{M}_I$, then $s(f)\in 1+\mathfrak{M}$, and $f^{-1}=\Phi(s(f))^{-1}$. 
\end{lemma}
\begin{proof}
Since $\langle f,1_{M/I}\rangle=1_R$, then $\langle s(f),1_M\rangle=1_R$ (because $1_{M/I}=1_M$). Therefore $s(f)\in 1+\mathfrak{M}$. Thus $s(f)^{-1}\in 1+\mathfrak{M}$, and $\Phi(s(f)^{-1})=\Phi(s(f))^{-1}$ (because $\Phi$ is a ring homomorphism). Finally, we obtain $f\Phi(s(f))^{-1}=\Phi(s(f))\Phi(s(f))^{-1}=1$ and $\Phi(s(f))^{-1}$ is a right inverse of $f$. The same holds for the left-side. 
\end{proof}

In the ordinary case, \textit{i.e.}, when $M$ is a (locally finite) monoid, the inverse $(-\zeta^{+})^*$ of the characteristic series $\zeta=1+\zeta^{+}$ is called the \emph{M\"obius series}, and denoted by $\mu(M)$. By analogy, we define the \emph{M\"obius series} of a locally finite monoid with zero $M$ as the series $\mu_0(M)=(-\zeta_0^{+})^*$, inverse of $\zeta_0=1+\zeta_0^{+}$ in $R_0[[M]]$. Therefore it satisfies $\mu_0(M)\zeta_0=\zeta_0\mu_0(M)=1$. 
\begin{lemma}\label{interlemma2}
Let $M$ be a locally finite monoid and $I$ be a two-sided proper ideal of $M$. Then, $\mu_0(M/I)=\Phi(\mu(M))$. Moreover if $\langle\mu(M),x\rangle=0_R$ for every $x \in I$, then $\mu_0(M/I)=\mu(M)$. 

\end{lemma}
\begin{proof}
The Rees quotient $M/I$ is a locally finite monoid with zero, and so its M\"obius series exists. Moreover $\zeta_0=\underline{M/I}_0 \in 1+\mathfrak{M}_I$, and according to Lemma~\ref{interlemma1}, $s(\zeta_0)\in 1+\mathfrak{M}$, and $(\zeta_0)^{-1}=\Phi(s(\zeta_0))^{-1}$. We have 
\begin{equation}\label{eqpourhilbertseries}
\begin{array}{lll}
s(\zeta_0)&=&\displaystyle\sum_{x\not\in I}x\\
&=&\underline{M\setminus I}\\
&=&\underline{M}-\underline{I}\\
&=&\zeta-\underline{I}\ .
\end{array}
\end{equation}
The series $\zeta^{+}-\underline{I}\in R[[M]]$ belongs to the augmentation ideal $\mathfrak{M}$ of $R[[M]]$ (as we already know), so the series $\zeta-\underline{I}=1+\zeta^{+}-\underline{I}$ is invertible in $R[[M]]$ with inverse $(\underline{I}-\zeta^{+})^*$. Therefore, according to Lemma~\ref{interlemma1}, 
\begin{equation}
\begin{array}{lll}
\mu_0(M/I)&=&\Phi(s(\zeta_0))^{-1}\\
&=&\Phi(s(\underline{M/I}_0))^{-1}\\
&=&\Phi(s(\underline{M/I}_0)^{-1})\\
&=&\Phi((\underline{I}-\underline{M^{+}})^*)\\
&=&\Phi((\underline{I}-\zeta^{+})^*)\\
&=&(\Phi(\underline{I}-\zeta^{+}))^* \\
&&\mbox{\begin{footnotesize}(because $\Phi$ is a continuous, for filtrations, algebra homomorphism)\end{footnotesize}}\\
&=&(\underbrace{\Phi(\underline{I})}_{=0}-\Phi(\zeta^{+}))^*\\
&=&\Phi((-\zeta^{+})^*)\\
&=&\Phi((1+\zeta^{+})^{-1})\\
&=&\Phi(\zeta^{-1})\\
&=&\Phi(\mu(M))\ .
\end{array}
\end{equation}
Now, if $\langle \mu(M),x\rangle=0_R$ for every $x \in I$, then $\mu_0(M)=\Phi(\mu(M))=\mu(M)$. 
\end{proof}

\begin{cor}\label{cor-pour-monoidelibre}
Let $X$ be any nonempty set. Let $I$ be a proper two-sided ideal of $X^*$. Then, 
\begin{equation}
\mu_0(X^*/I)=\left \{
\begin{array}{lll}
\mu(X^*) & \mbox{if} & X\cap I=\emptyset\ ,\\
\mu((X\setminus I)^*) & \mbox{if} &X\cap I\not=\emptyset\ .
\end{array}
\right .
\end{equation}
\end{cor}
\begin{proof}
We can apply Lemma~\ref{interlemma2} to obtain $\mu_0(X^*/I)=\Phi(\mu(X^*))$. According to~\cite{cartfoa}, $\mu(X^*)=1-\sum_{x\in X}x$. If $X \cap I=\emptyset$, then $\mu_0(X^*/I)=\Phi(\mu(X^*))=\mu(X^*)$, and if $X\cap I \not=\emptyset$, then let $Y=X\setminus I$. We have $\Phi(\mu(X^*))=1-\sum_{y\in Y}y=\mu(Y^*)$. 
\end{proof}

\begin{example}
\begin{enumerate}
\item Let $X$ be any nonempty set. Let $I=\{\omega \in X^* : \exists x\in X,\ |\omega|_x\geq 2\}$. The set $X^*/I$ consists of all standard words, \textit{i.e.}, word without repetition of any letter. Then according to Corollary~\ref{cor-pour-monoidelibre}, $\mu_0(X^*/I)=\mu(X^*)=1-\underline{X}$ as announced in Section~\ref{section-introduction} Introduction. 
\item Let $X$ be any set. A congruence $\equiv$ on $X^*$ is said to be \emph{multihomogeneous}~\cite{DK,DK94} if, and only if, $\omega\equiv \omega'$ implies $|\omega|_x=|\omega'|_x$ for every $x \in X$. A quotient monoid $X^*/\equiv$ of $X^*$ by a multihomogeneous congruence is called a \emph{multihomogeneous monoid}. For instance, any free partially commutative monoid, the plactic~\cite{LS81}, hypoplactic~\cite{KT97,Nov00}, Chinese~\cite{DK94} and sylvester~\cite{HNT02} monoids are multihomogeneous. Such a monoid $M=X/\equiv$ is locally finite and therefore admits a M\"obius function $\mu$ with $\mu(1_M)=1$ and $\mu(x)=-1$ for every $x \in X$. An epimorphism $\mathsf{Ev}$ from $M$ onto the free commutative monoid $X^{\oplus}$, the \emph{commutative image}, is given by $\mathsf{Ev}(\omega)=\sum_{x\in X}|\omega|_x \delta_x$, where $\delta_x$ is the indicator function of $x$. Any proper ideal $I$ of $X^{\oplus}$ gives rise to a proper two-sided ideal $\mathsf{Ev}^{-1}(I)$ of $M$. Let $I=\{f\in X^{\oplus} : \sum_{x\in X}f(x)\geq 2\}$. Then, as sets, $M/\mathsf{Ev}^{-1}(I)=\{0,1_M\}\cup X$ and  $\mu_0(M/\mathsf{Ev}^{-1}(I))=1-\underline{X}$. 
\end{enumerate}
\end{example}

\section{Some remarks about Hilbert series}
Now, let $X$ be a finite set, and $I$ be a proper two-sided ideal of $X^*$. For any $S\subseteq X^*$ or $S\subseteq X^*/ I$, we define $S(n)=\{w \in S : |w|=n\}$ for any $n\in\mathbb{N}$. (Note that the notation $M(n)$ is consistent with the given one in remark~\ref{equivalencetopologies} for $M=X^*/I$.) Let $\mathbb{K}$ be a field. Let define $A_n=\mathbb{K} X^*(n)$ (the $\mathbb{K}$-vector space spanned by $X^*(n)$), and $B_n=\mathbb{K} (X^*/I)(n)$ for every $n\in \mathbb{N}$, in such a way that $\mathbb{K}[X^*]=\bigoplus_{n\geq 0} A_n$ and $\mathbb{K}_0[X^*/ I]=\bigoplus_{n\geq 0}B_n$. (Note that for every $w,w'\in X^*/I$, we have $|ww'|=|w|+|w'|$ if $ww'\not=0$, in such a way that $B_m B_n\subseteq B_{m+n}$ since $0\in B_i$ for every $i$.) Since $X$ is finite, for every integer $n$, $X^*(n)$ and $(X^*/ I)(n)$ are finite, and therefore 
$A_n$ and $B_n$ are finite-dimensional $\mathbb{K}$-vector spaces. Moreover $\dim(B_n)=\dim(A_n)-|I(n)|$ because $(X^*/I)(n)=X^*(n)\setminus I(n)$. So in particular respective Hilbert series\footnote{Let $\mathcal{A}=\bigoplus_{n\geq 0}A_n$ be a graded algebra, where for every $n$, $A_n$ is finite dimensional. The \emph{Hilbert series of $\mathcal{A}$} (in the variable $t$) is defined by $\mathpzc{Hilb}_{\mathcal{A}}(t)=\sum_{n\geq 0}\dim(A_n)t^n$.} are related by 
\begin{equation}
\displaystyle \mathpzc{Hilb}_{\mathbb{K}_0[X^*/ I]}(t)=\mathpzc{Hilb}_{\mathbb{K}[X]}(t)-\sum_{n\geq 0}|I(n)|t^n=\frac{1}{1-|X|t}-\sum_{n\geq 0}|I(n)|t^n\ .
\end{equation}
Note that since $I$ is a proper ideal, $I(0)=\emptyset$, and $\sum_{n\geq 1}|I(n)|t^n$ may be interpreted rather naturally as the Hilbert series of the ideal $\mathbb{K}[I]=\bigoplus_{n\geq 1}\mathbb{K}I(n)$ (it also follows that $\mathpzc{Hilb}_{\mathbb{K}[I]}(t)$ is not invertible in $\mathbb{Z}[[t]]$). We have 
\begin{equation}
\mathpzc{Hilb}_{\mathbb{K}_0[X^*/ I]}(t)=\mathpzc{Hilb}_{\mathbb{K}[X]}(t)-\mathpzc{Hilb}_{\mathbb{K}[I]}(t)\ .
\end{equation}
This equation may be recovered from equation~(\ref{eqpourhilbertseries}), namely $s(\zeta_0)=\zeta-\underline{I}$, using an evaluation map. Suppose now that $\mathbb{K}$ is a field of characteristic zero, and $t$ is a variable. Let $e:X^*\rightarrow \{t^i\}_{i\in\mathbb{N}}$ be the unique morphism of monoids such that $e(x)=t$ for every $x\in X$. We extend it to a $\mathbb{Z}$-linear map from $\mathbb{Z}[X^*]$ to $\mathbb{Z}[t]$ by $e(\sum_{w\in X^*}n_w w)=\sum_{n\in \mathbb{N}}\left( \sum_{w\in X^*,\ |w|=n}n_w\right )t^n$ from . Moreover since $X$ is finite, for every $n\in\mathbb{N}$, $X(n)$ is also finite (of cardinality $|X|^n$), and therefore for every series $f=\sum_{w\in X^*}n_w w\in \mathbb{Z}[[X^*]]$, by summability, we have $f=\sum_{n\in \mathbb{N}}f_n$, where $f_n=\sum_{w\in X^*(n)}n_w w\in \mathbb{Z}[X^*]$ for every $n\in \mathbb{N}$, and we extend $e$ (by continuity) as a linear map from $\mathbb{Z}[[X^*]]$ to $\mathbb{Z}[[t]]$ by $e(f)=\sum_{n\in \mathbb{N}}e(f_n)=\sum_{n\in\mathbb{N}}\left (\sum_{w\in X^*(n)}n_w \right )t^n$. Now, applying $e$ on both side of equation~(\ref{eqpourhilbertseries}), we obtain (note that $s(\zeta_0), \zeta$ and $\underline{I}$ belong to $\mathbb{Z}[[X^*]]$) 
\begin{equation}
\begin{array}{llll}
&e(s(\zeta_0))&=&e(\zeta)-e(\underline{I})\\
\Leftrightarrow & e\left (\sum_{w\not\in I}w\right )&=&e\left (\sum_{w\in X^*}w\right)-e\left(\sum_{w\in I}w\right )\\
\Leftrightarrow & e\left(\sum_{n\in \mathbb{N}}\underline{X(n)\setminus I(n)}\right )&=&e\left ( \sum_{n\in\mathbb{N}}\underline{X(n)}\right )-e\left (\sum_{n\in\mathbb{N}}\underline{I(n)}\right)\\
\Leftrightarrow & \sum_{n\in\mathbb{N}}(|X(n)|-|I(n)|)t^n&=&\sum_{n\in\mathbb{N}}|X(n)|t^n-\sum_{n\in\mathbb{N}}|I(n)|t^n\\
\Leftrightarrow & 
\mathpzc{Hilb}_{\mathbb{K}_0[X^*/ I]}(t)&=&\mathpzc{Hilb}_{\mathbb{K}[X]}(t)-\mathpzc{Hilb}_{\mathbb{K}[I]}(t)\ .
\end{array}
\end{equation}
Last equality is nothing else than the obvious relation between the ordinary generating functions of the combinatorial class $X^*\setminus I$, $X^*$ and $I$, where the notion of size is nothing else than the length of words (see~\cite{Fla09}, theorem I.5 `` implicit specifications ''). 
\begin{example}
\begin{enumerate}
\item Suppose that $I=\{\omega \in X^* : \exists x\in X,\ |\omega|_x\geq 2\}$. It is clear that for every $n> |X|$, $I(n)=X(n)$. For every $n\leq|X|$, $|(X^*/ I)(n)|=\prod_{i=0}^{n-1}(|X|-i)=|X|^{\underline{n}}$ (in particular,  $|(X^*/ I)(0)|=|\{\epsilon\}|=1$, and $|(X^*/ I)(1)|=|X|$). If follows that $\mathpzc{Hilb}_{\mathbb{K}_0[X^*/ I]}(t)=\sum_{n=0}^{|X|}|X|^{\underline{n}}t^n$, and therefore 
$\mathpzc{Hilb}_{\mathbb{K}[I]}(t)=\sum_{n\geq 2}(|X|^n-|X|^{\underline{n}})t^n$. 
\item Let $n_0\in\mathbb{N}$ such that $n_0\geq 1$. Let $I=\{w\in X^* : |w|>n_0\}$. Then $\mathpzc{Hilb}_{\mathbb{K}_0[X^*/I]}(t)=\sum_{n=0}^{n_0}|X|^nt^n$, $\mathpzc{Hilb}_{\mathbb{K}[I]}(t)=\sum_{n\geq n_0+1}|X|^nt^n$. 
\end{enumerate}
\end{example}

\bibliographystyle{plain}
\bibliography{biblio}








\end{document}